%% file: AnEfficientMethodOfSplineApproximation.tex
\documentclass[12pt,letterpaper,oneside,reqno]{amsart}
\usepackage{amsfonts}
\usepackage{amsmath}
\usepackage{amssymb}
\usepackage{amsthm}
\usepackage{float}
\usepackage{mathrsfs}
\usepackage{colonequals}
\usepackage[font=small,labelfont=bf]{caption}
\usepackage[left=1in,right=1in,bottom=1in,top=1in]{geometry}
\usepackage[pdfpagelabels,hyperindex,colorlinks=true,linkcolor=blue,urlcolor=magenta,citecolor=green]{hyperref}
\usepackage{graphicx}
\usepackage{array}
\usepackage{etoolbox}
\apptocmd{\sloppy}{\hbadness 10000\relax}{}{}
\newcommand \coeffA [3][A] {{\mathbf{#1}} \sb{#2,#3}}

\let\svthefootnote\thefootnote
\newcommand\freefootnote[1]{%
    \let\thefootnote\relax%
    \footnotetext{#1}%
    \let\thefootnote\svthefootnote%
}

\title[An efficient method of spline approximation for power function]
{An efficient method of spline approximation for power function}
\author[Petro Kolosov]{Petro Kolosov}
\address{Software Developer, DevOps Engineer}
\email{kolosovp94@gmail.com}
\urladdr{https://kolosovpetro.github.io}
\date{\today}
\hypersetup{
    pdftitle={An efficient method of spline approximation for power function},
    pdfsubject={
        Polynomials,
        Finite differences,
        Interpolation,
        Approximation,
        Polynomial identities,
        Power sums,
        Binomial theorem,
        Power function,
        Binomial coefficients,
        Bernoulli numbers,
        Pascal's triangle,
        Faulhaber's formula,
        OEIS,
        Bernoulli polynomials,
        Combinatorics,
        Approximation,
        Splines,
        Spline approximation
    },
    pdfauthor={Petro Kolosov},
    pdfkeywords={
        Polynomials,
        Finite differences,
        Interpolation,
        Approximation,
        Polynomial identities,
        Power sums,
        Binomial theorem,
        Power function,
        Binomial coefficients,
        Bernoulli numbers,
        Pascal's triangle,
        Faulhaber's formula,
        OEIS,
        Bernoulli polynomials,
        Combinatorics,
        Approximation,
        Splines,
        Spline approximation
    }
}
\begin{document}
    \begin{abstract}
        \input{sections/01_abstract}
    \end{abstract}

    \maketitle

    \tableofcontents

    \freefootnote{Sources: \url{https://github.com/kolosovpetro/AnEfficientMethodOfSplineApproximation}}

    \section{Introduction} \label{sec:introduction}
    \input{sections/02_introduction}

    \section{Generalizations}\label{sec:generalizations}
    \input{sections/03_generalizations}

    \section{Use cases}\label{sec:use-cases}
    \input{sections/04_use_cases}

    \section{Conclusions}\label{sec:conclusions}
    \input{sections/conclusions}

    \bibliographystyle{unsrt}
    \bibliography{AnEfficientMethodOfSplineApproximation}
    \noindent \textbf{Version:} \input{sections/version}
\end{document}

%% file: sections/01_abstract.tex
Let $P(m, X, N)$ be an $m$-degree polynomial in $X\in\mathbb{R}$
having fixed non-negative integers $m$ and $N$.
The polynomial $P(m, X, N)$ is derived from a rearrangement of Faulhaber's formula
in the context of Knuth's work entitled "Johann Faulhaber and sums of powers".
In this manuscript we discuss the approximation properties of polynomial $P(m,X,N)$.
In particular, the polynomial $P(m,X,N)$ approximates the odd power function $X^{2m+1}$ in a certain neighborhood
of a fixed non-negative integer $N$ with a percentage error under $1\%$.
By increasing the value of $N$ the length of convergence interval with odd-power $X^{2m+1}$ also increases.
Furthermore, this approximation technique is generalized for arbitrary non-negative exponent $j$ of the power function $X^j$
by using splines.

%% file: sections/02_introduction.tex
Consider the $m$-degree polynomial $P(m, X, N)$ having fixed non-negative integers $m$ and $N$
\begin{align*}
    P(m,X,N) = \sum_{r=0}^{m} \sum_{k=1}^{N} \coeffA{m}{r} k^r (X-k)^r
\end{align*}
For example
\begin{align*}
    P(2,X,0) &= 0 \\
    P(2,X,1) &= 30X^2 - 60X + 31 \\
    P(2,X,2) &= 150X^2 - 540X + 512 \\
    P(2,X,3) &= 420X^2 - 2160X + 2943 \\
    P(2,X,4) &= 900X^2 - 6000X + 10624
\end{align*}
where $\coeffA{m}{r}$ is a real coefficient defined recursively, see~\cite{alekseyev2018mathoverflow,
    on_the_link_between_binomial_theorem_and_discrete_convolution, unusual_identity_for_odd_powers,
    history_and_overview_of_polynomial_p}.
For example,
\input{sections/figures/05_fig_coefficients_a}

The polynomial $P(m, X, N)$ is derived from a rearrangement of Faulhaber's formula
in the context of Knuth's work entitled \textit{Johann Faulhaber and sums of powers}, see~\cite{knuth1993johann}.
In particular, the polynomial $P(m, X, N)$ yields an identity for odd powers
\begin{align*}
    P(m, X, X) = X^{2m+1}
\end{align*}
In its extended form
\begin{align*}
    X^{2m+1} = \sum_{r=0}^{m} \sum_{k=1}^{X} \coeffA{m}{r} k^r (X-k)^r
\end{align*}
The exact relation between Faulhaber's formula and $P(m,X,N)$ is shown by~\cite{kolosov2025unexpected}.

However, apart from the polynomial identity for odd powers, a few approximation properties of $P(m,X,N)$
were discovered in addition.
Therefore, in this manuscript we explore the approximation properties of the polynomial $P(m,X,N)$.
During our discussion, we utilize the following well-known criteria to measure and estimate
the error of approximation: Absolute error, Relative error and Percentage error.
Assume that the function $f_2(x)$ approximates the function $f_1 (x)$, then errors are given by
\begin{align*}
    \mathrm{Absolute \; Error}   &= \lvert f_1(x) - f_2(x) \rvert \\
    \mathrm{Relative \; Error}   &= \frac{\lvert f_1(x) - f_2(x) \rvert}{\lvert f_1(x) \rvert} \\
    \mathrm{Percentage \; Error} &= \frac{\lvert f_1(x) - f_2(x) \rvert}{\lvert f_1(x) \rvert} \times 100\%
\end{align*}

Diving straight into the point, we switch our focus to the partial case of polynomial
$P(2,X,4) = 900X^2 - 6000X + 10624$
to show the first example of how it approximates the odd power function $X^5$.
In general, we approximate the polynomial $X^{2m+1}$ using a lower-degree polynomial of degree $m$.
The following image demonstrates the approximation of fifth power $X^5$ by
$P(2,X,4) = 900X^2 - 6000X + 10624$
\begin{figure}[H]
    \centering
    \includegraphics[width=1\textwidth]{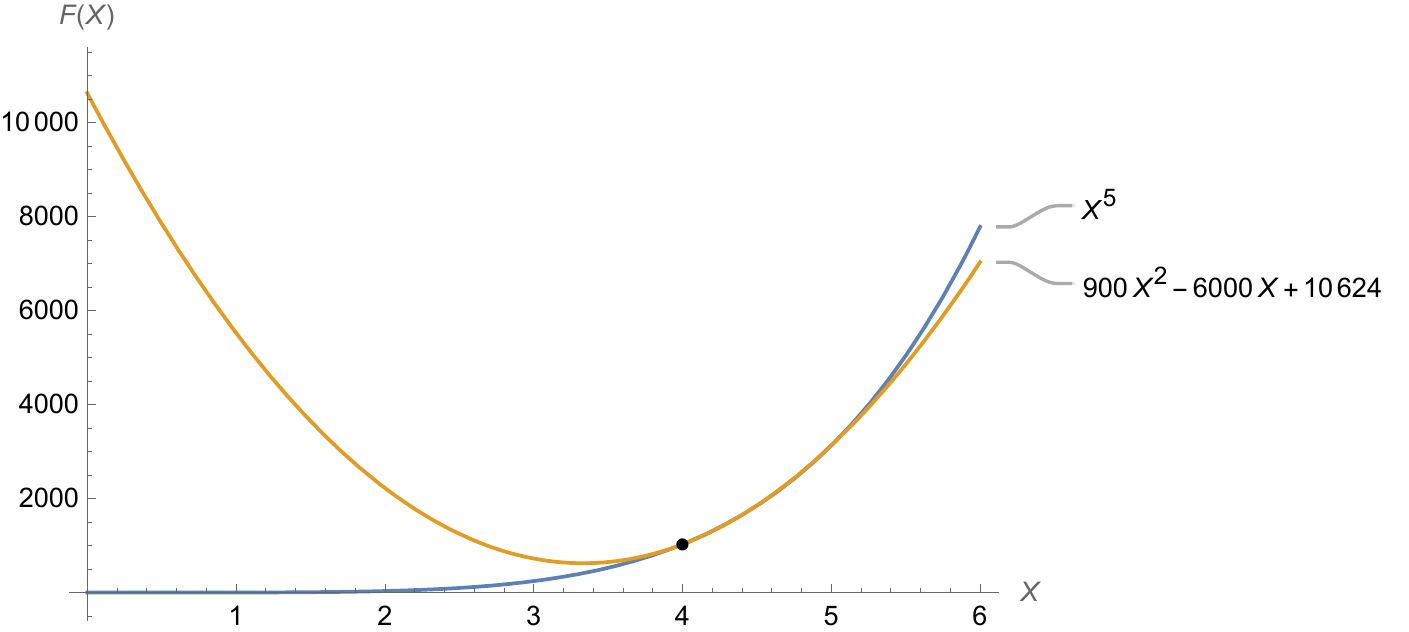}
    ~\caption{Approximation of fifth power $X^5$ by $P(2, X, 4)$.
    Convergence interval is $4.0 \leq X \leq 5.1$ with a percentage error $E < 1\%$.
    }\label{fig:03_plots_polynomial_p2_n4_with_fifth}
\end{figure}
As observed, the polynomial $P(2, X, 4)$ approximates $X^5$ in a certain neighborhood of $N=4$ with
the convergence interval $4.0 \leq X \leq 5.1$ such that the percentage error is less than $1\%$ which is quite remarkable.
The following table presents specific values of absolute, relative, and percentage errors for this approximation
\input{sections/figures/032_polynomials_p2_table_n4}

One more interesting observation arises by increasing the value of $N$ in $P(m, X, N)$ while keeping $m$ fixed.
As $N$ increases, the length of the convergence interval with the odd-power $X^{2m+1}$ also increases.
For instance,
\begin{itemize}
    \item For $P(2, X, 4)$ and $X^5$, the convergence interval with a percentage error less than $1\%$ is $4.0 \leq X \leq 5.1$, with a length $L=1.1$
    \item For $P(2, X, 20)$ and $X^5$, the convergence interval with a percentage error less than $1\%$ is $18.7 \leq X \leq 22.9$, with a length $L=4.2$
    \item For $P(2, X, 120)$ and $X^5$, the convergence interval with a percentage error less than $1\%$ is $110.0 \leq X \leq 134.7$, with a length $L=24.7$
\end{itemize}

The reason behind this behavior lies in the implicit form of the polynomial $P(m,X,N)$,
meaning that
\begin{align*}
    P(m,X,N) = \sum_{r=0}^{m} (-1)^{m-r} U(m, N, r) \cdot X^{r}
\end{align*}
where $U(m, N, r)$ is a polynomial defined as follows
\begin{align*}
    U(m, N, r) = (-1)^m \sum_{k=1}^{N} \sum_{j=r}^{m} \binom{j}{r} \coeffA{m}{j} k^{2j-r} (-1)^j
\end{align*}
which grows as $N$ increases.
Few cases of coefficients $U(m, N, r)$ are registered as OEIS sequences~\cite{
    oeis_coefficients_u_m_l_k_defined_by_polynomial_identity_1,
    oeis_coefficients_u_m_l_k_defined_by_polynomial_identity_2,
    oeis_coefficients_u_m_l_k_defined_by_polynomial_identity_3}.

To summarize, let us recap the key findings so far.
The polynomial $P(m,X,N)$ is an $m$-degree polynomial in $X \in \mathbb{R}$ with fixed non-negative integers $m$ and $N$.
It approximates the odd power function $X^{2m+1}$ within a specific neighborhood of $N$.
The length $L$ of the convergence interval between $X^{2m+1}$ and $P(m,X,N)$ increases as $N$ grows.

For the sake of clear and precise verification of results, consider the Mathematica programs to generate
plots and data tables, so that reader is able to verify the main results of current part of manuscript,
see~\cite{kolosovpetro_gist}.

So far we have discussed approximation of odd power function $X^{2m+1}$, now we focus on its even case $X^{2m+2}$
which is quite straightforward.
Considering the same example $P(2, X, 4)$ we reach the approximation of even power $X^6$
by means of $K$-times multiplication by $X$, with graphic representation as follows
\begin{figure}[H]
    \centering
    \includegraphics[width=1\textwidth]{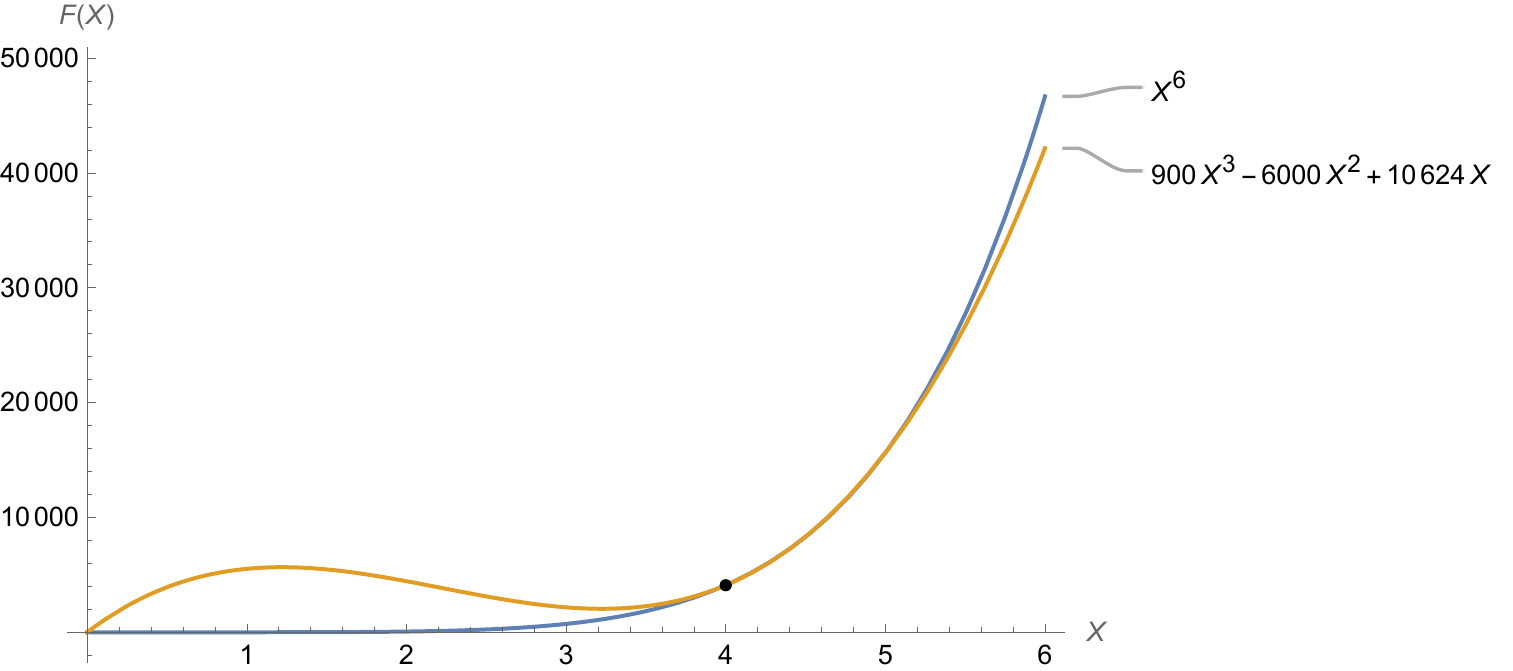}
    ~\caption{Approximation of sixth power $X^6$ by $P(2, X, 4) \cdot X$.
    Convergence interval is $3.9 \leq X \leq 5.1$ with a percentage error $E < 3\%$.
    }\label{fig:07_plot_of_6th_power_with_p_2_4_times_x}
\end{figure}
Therefore, we have reached the statement that
the polynomial $P(m,X,N)$ is an $m$-degree polynomial in $X$, having fixed non-negative
integers $m$ and $N$.
It approximates the power function $X^{j}$ in a certain neighborhood of fixed $N$.
The length of convergence interval between the power function $X^j$ and $P(m,X,N) \cdot X^K$ increases as $N$ grow.

%% file: sections/figures/05_fig_coefficients_a.tex
\begin{table}[H]
    \begin{center}
        \setlength\extrarowheight{-6pt}
        \begin{tabular}{c|cccccccc}
            $m/r$ & 0 & 1       & 2      & 3      & 4   & 5    & 6     & 7 \\ [3px]
            \hline
            0     & 1 &         &        &        &     &      &       &       \\
            1     & 1 & 6       &        &        &     &      &       &       \\
            2     & 1 & 0       & 30     &        &     &      &       &       \\
            3     & 1 & -14     & 0      & 140    &     &      &       &       \\
            4     & 1 & -120    & 0      & 0      & 630 &      &       &       \\
            5     & 1 & -1386   & 660    & 0      & 0   & 2772 &       &       \\
            6     & 1 & -21840  & 18018  & 0      & 0   & 0    & 12012 &       \\
            7     & 1 & -450054 & 491400 & -60060 & 0   & 0    & 0     & 51480
        \end{tabular}
    \end{center}
    \caption{Coefficients $\coeffA{m}{r}$. See the OEIS sequences
    ~\cite{oeis_numerators_of_the_coefficient_a_m_r,oeis_denominators_of_the_coefficient_a_m_r}.}
    \label{tab:table_of_coefficients_a}
\end{table}

%% file: sections/figures/032_polynomials_p2_table_n4.tex
\begin{table}[H]
    \centering
    \begin{tabular}{|c|c|c|c|c|c|}
        \hline
        \textbf{X} & \textbf{$X^5$} & \textbf{$900X^2 - 6000X + 10624$} & \textbf{ABS} & \textbf{Relative} & \textbf{\% Error} \\ \hline
        4.0        & 1024.0         & 1024.0                            & 0.0          & 0.0               & 0.0               \\ \hline
        4.1        & 1158.56        & 1153.0                            & 5.56201      & 0.00480079        & 0.480079          \\ \hline
        4.2        & 1306.91        & 1300.0                            & 6.91232      & 0.00528905        & 0.528905          \\ \hline
        4.3        & 1470.08        & 1465.0                            & 5.08443      & 0.0034586         & 0.34586           \\ \hline
        4.4        & 1649.16        & 1648.0                            & 1.16224      & 0.000704746       & 0.0704746         \\ \hline
        4.5        & 1845.28        & 1849.0                            & 3.71875      & 0.00201528        & 0.201528          \\ \hline
        4.6        & 2059.63        & 2068.0                            & 8.37024      & 0.00406395        & 0.406395          \\ \hline
        4.7        & 2293.45        & 2305.0                            & 11.5499      & 0.00503605        & 0.503605          \\ \hline
        4.8        & 2548.04        & 2560.0                            & 11.9603      & 0.00469393        & 0.469393          \\ \hline
        4.9        & 2824.75        & 2833.0                            & 8.24751      & 0.00291973        & 0.291973          \\ \hline
        5.0        & 3125.0         & 3124.0                            & 1.0          & 0.00032           & 0.032             \\ \hline
        5.1        & 3450.25        & 3433.0                            & 17.2525      & 0.00500036        & 0.500036          \\ \hline
    \end{tabular}
    \caption{Comparison of $X^5$ and $P(2,X,4) = 900X^2 - 6000X + 10624$}
    \label{tab:table2}
\end{table}

%% file: sections/03_generalizations.tex
Previously, we have discussed that polynomial $P(m,X,N)$ approximates power function $X^j$
in some neighborhood of fixed non-negative integer $N$.
Approximation by $P(m,X,N)$ can be adjusted by $X^k$ multiplication for even exponents of power function.
In general, it is safe to say that power function $X^j$ is approximated by $P(m,X,N) \cdot X^k$
where $k=0$ for odd exponent $j$ and $k$ is either $k=1$ or $k=-1$ for an even exponent $j$.
Therefore, for arbitrary exponent $j$ in $X^j$ we have
\begin{align*}
    X^j \approx
    \begin{cases}
        P(m,X,N) \quad         & j=2m+1 \\
        P(m,X,N) \cdot X \quad & j=2m+2 \\
        P(m,X,N) \cdot X^{-1} \quad & j=2m \\
    \end{cases}
\end{align*}
Of course, there are other variations of the value of $k$, but we will stick to the simple case for the moment.

As we also discussed, the length $L$ of the convergence interval between $X^j$ and its approximation by $P(m,X,N)$
increases as $N$ grow.
However, the convergence interval is still bounded, which could not satisfy certain approximation scenarios.
Depending on the approximation requirements in terms of convergence interval length $L$ a single polynomial $P(m,X,N)$
with fixed $m$ and $N$ may be unsuitable.
Here is the place where spline approximation comes into play.
The spline $S(x)$ is piecewise defined function over the interval $(x_0, \ldots x_n)$
\begin{align*}
    S(x) &=
    \begin{cases}
        f_1(x), & x_0 \leq x < x_1 \\
        f_2(x), & x_0 \leq x < x_1\\
        \vdots & \vdots \\
        f_n(x), & x_{n-1} \leq x \leq x_n
    \end{cases}
\end{align*}
The given points $x_k$ are called \textit{knots}.

Assume that the approximation requirement in terms of convergence length $L$ is to approximate the power function $X^j$
bounded by real points $A$ and $B$ such that $A < B$.
Splines perfectly fit the need to match an arbitrary convergence range for the power function $X^j$ using the
approximation by $P(m,X,N)$.
Formally,
\begin{align*}
    X^j \approx
    \begin{cases}
        P(m,X,N+t_1) \cdot X^{k} \quad & x_0 \leq x < x_1 \\
        P(m,X,N+t_2) \cdot X^{k} \quad & x_0 \leq x < x_1 \\
        \vdots & \vdots \\
        P(m,X,N+t_{n-1}) \cdot X^{k} \quad & x_{n-1} \leq x < x_n
    \end{cases}
\end{align*}
The values of $t_r$ to be adjusted according to approximation requirements in terms of accuracy.

%% file: sections/04_use_cases.tex
The approximation technique above has its own constraints and limitations.
For instance, approximation requirements should have precisely specified exponent $j$ in $X^j$ because
for each $j$ there is a matching polynomial $P(m,X,N)$.
Perfectly, there should be a set precompiled polynomials $P(m,X,N)$ matching precise exponent $j$ in $X^j$ over
exactly defined approximation range with required error of approximation $E$ as a constraint.
Generally, the approximation of power function $X^j$ by $P(m,X,N)$ can be broken down into the following steps
\begin{enumerate}
    \item Define the exponent $j$ in $X^j$
    \item Define the required error threshold $E$
    \item Define the required interval of approximation $I$
    \item Choose and precompile polynomials $P(m,X,N)$
    so that the required interval of approximation and error threshold $E$ are satisfied
    \item Define a set of knots so that the error threshold $E$ and the interval of approximation $I$ requirements are satisfied
\end{enumerate}
Defining set of spline knots essentially requires an inspection of the
convergence intervals between $X^j$ and $P(m,X,N+t_k)$
by choosing knots such that the interval of approximation and error threshold are satisfied.
Consider an example.
Let be the following approximation requirements
\begin{enumerate}
    \item Exponent $j=3$
    \item Percentage error threshold $E\leq 1\%$
    \item Interval of approximation $10 \leq X \leq 15$
\end{enumerate}
Now we have to choose a set of polynomials $P(m, X, N+t_k)$ based on which we adjust a set of spline knots.
We can safely choose integers $t_k$ in range $10 \leq t_k \leq 15$ because
of the following properties of $P(m,X, N)$.
\begin{align*}
    P(m,X, X) &= X^{2m+1} \\
    P(m,X, X+1) &= (X+1)^{2m+1} - 1
\end{align*}
Therefore, for each two consequential points $N=X, N=X+1$ the absolute difference $P(m,X, X) - P(m,X, X+1)$ equals to one,
consequently, the range $[X, X+1]$ satisfies the $1\%$ error threshold for any pair $j, X$ such that $X^j \leq 100$.
The approximation range $10 \leq X \leq 15$ and exponent $j=3$ are chosen intentionally to show the spline approximation with
percentage error threshold less than $1\%$.
Therefore, to approximate the polynomial $X^3$ in the range $10 \leq X \leq 15$, we use the following spline function
\begin{align}
    X^3 \approx
    \begin{cases}
        P(1,X,10) = -2300 + 330X, & 10 \leq X < 11 \\
        P(1,X,11) = -3025 + 396X, & 11 \leq X < 12 \\
        P(1,X,12) = -3888 + 468X, & 12 \leq X < 13 \\
        P(1,X,13) = -4901 + 546X, & 13 \leq X < 14 \\
        P(1,X,14) = -6076 + 630X, & 14 \leq X \leq 15
    \end{cases}
    \label{eq:spline_approximation_of_cubes}
\end{align}
Graphically, this linear approximation of cubes looks as follows
\begin{figure}[H]
    \centering
    \includegraphics[width=1\textwidth]{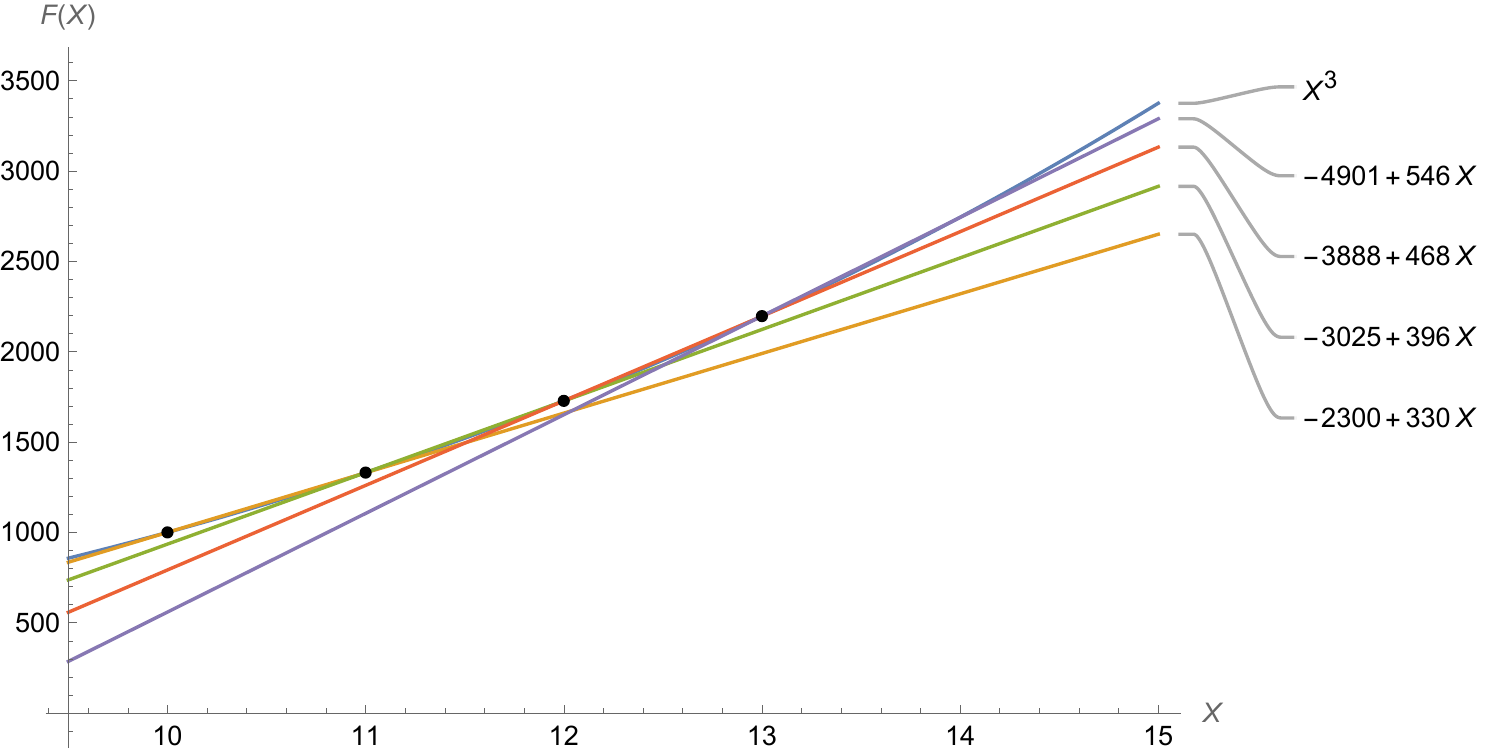}
    ~\caption{
        Approximation of cubes $X^3$ by splines~\eqref{eq:spline_approximation_of_cubes}.
        Convergence interval is $10 \leq X \leq 15$ with a percentage error $E < 1\%$.
    }
    \label{fig:08_plots_of_cubes_power_with_p_2_10_15}
\end{figure}
where the spline knots are integers in the range $10 \leq N \leq 14$.

The same principle applies for even exponent $j=4$ in $X^j$ with the same convergence interval $10 \leq X \leq 15$
and approximation error under $1\%$
\begin{align}
    X^4 \approx
    \begin{cases}
        P(1,X,10) \cdot X = -2300X + 330X^2, & 10 \leq X < 11 \\
        P(1,X,11) \cdot X = -3025X + 396X^2, & 11 \leq X < 12 \\
        P(1,X,12) \cdot X = -3888X + 468X^2, & 12 \leq X < 13 \\
        P(1,X,13) \cdot X = -4901X + 546X^2, & 13 \leq X < 14 \\
        P(1,X,14) \cdot X = -6076X + 630X^2, & 14 \leq X \leq 15
    \end{cases}
    \label{eq:spline_approximation_fourth_power}
\end{align}
Which graphically looks as follows
\begin{figure}[H]
    \centering
    \includegraphics[width=1\textwidth]{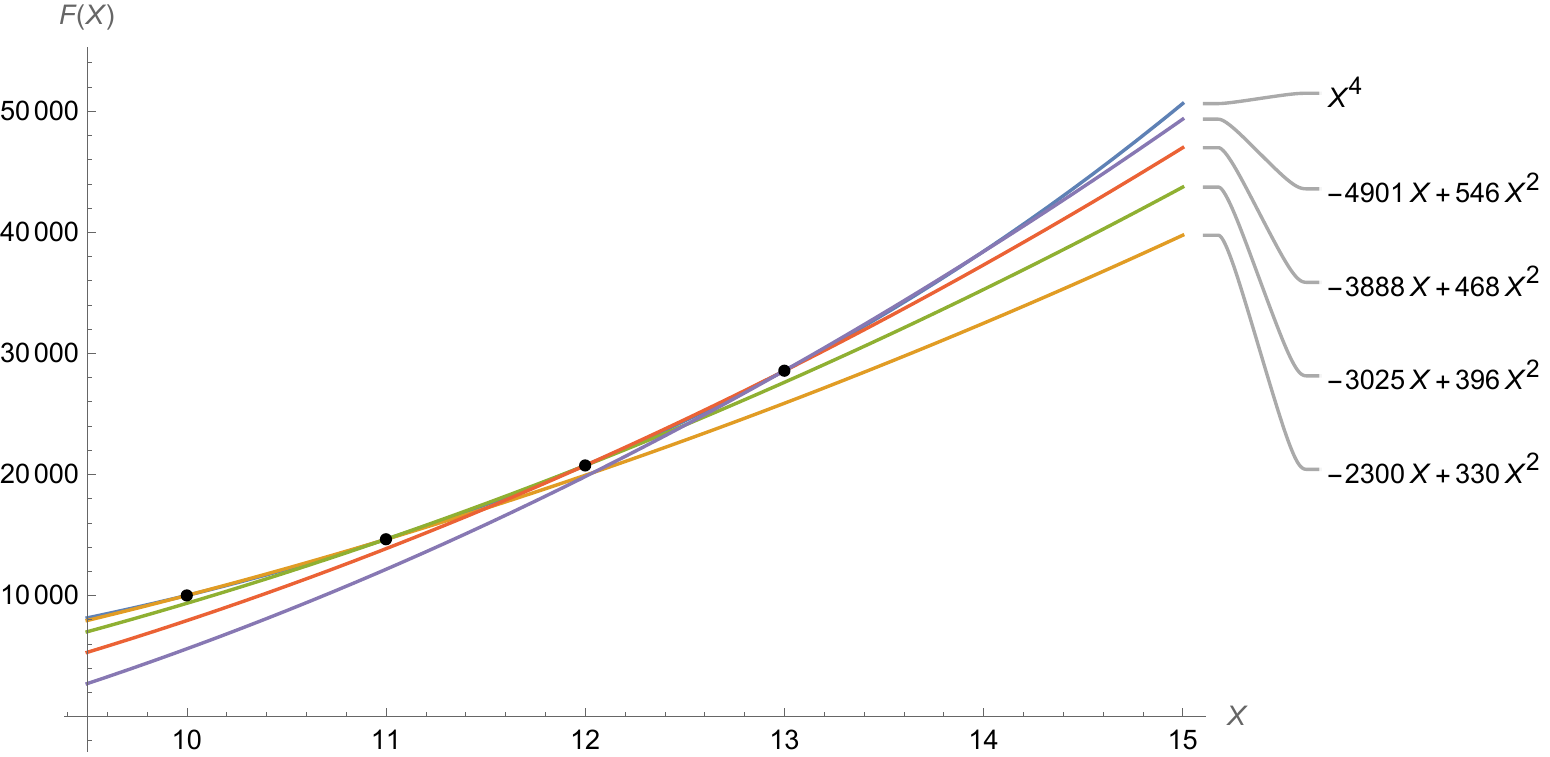}
    ~\caption{
        Approximation of $X^4$ by splines~\eqref{eq:spline_approximation_fourth_power}.
        Convergence interval is $10 \leq X \leq 15$ with a percentage error $E < 1\%$.
    }
    \label{fig:09_plots_of_fourth_power_with_p_2_10_15_times_x}
\end{figure}
In general, for each variation of $X^j$ such that $X^j \geq 100$ the approximation can be done using
splines in $P(m,X, N)$ over the interval $A \leq X \leq B$ with spline knot vector be the integers in
range from $A$ to $B$ so that knots vector is $\{A, A+1, A+2, \ldots, B \}$.
Because,
\begin{align*}
    P(m,X, X) &= X^{2m+1} \\
    P(m,X, X+1) &= (X+1)^{2m+1} - 1
\end{align*}
This can be further optimized depending on the value of $N$ in $P(m,X,N)$ because the convergence interval
with the power function $X^j$ increases as $N$ grows.

%% file: sections/conclusions.tex
We have established that $P(m, X, N)$ is an $m$-degree polynomial in $X\in\mathbb{R}$
having fixed non-negative integers $m$ and $N$.

The polynomial $P(m, X, N)$ is a result of a rearrangement in Faulhaber's formula
in the context of Knuth's work \textit{Johann Faulhaber and sums of powers}, see~\cite{knuth1993johann}.
The exact relation between Faulhaber's formula and $P(m,X,N)$ is shown by~\cite{kolosov2025unexpected}.

In this manuscript we have discussed the approximation properties of the polynomial $P(m,X,N)$.

In particular, the polynomial $P(m,X,N)$ approximates odd power function $X^{2m+1}$ in a certain neighborhood
around a fixed non-negative integer $N$ with percentage error of less than $1\%$.

By increasing the value of $N$ the length of the convergence interval with the odd-power $X^{2m+1}$ also increases.

Furthermore, this approximation technique is generalized to an arbitrary non-negative exponent power function $X^j$
by using splines.

In general, for each variation of $X^j$ such that $X^j \geq 100$ the approximation can be done using
splines in $P(m,X, N)$ over the interval $A \leq X \leq B$ with the spline knot vector being the integers in
range from $A$ to $B$, so that the knots vector is $\{A, A+1, A+2, \ldots, B \}$.
Because,
\begin{align*}
    P(m,X, X) &= X^{2m+1} \\
    P(m,X, X+1) &= (X+1)^{2m+1} - 1
\end{align*}
This can be further optimized depending on the value of $N$ in $P(m,X,N)$ because the convergence interval
with the power function $X^j$ increases as $N$ grows.

%% file: sections/version.tex
\texttt{1.0.5-tags-v1-0-4.5+tags/v1.0.4.387ee0e}

%% file: AnEfficientMethodOfSplineApproximation.bbl
\begin{thebibliography}{10}

\bibitem{alekseyev2018mathoverflow}
{Alekseyev, Max}.
\newblock {MathOverflow answer 297916/113033}, 2018.
\newblock \url{https://mathoverflow.net/a/297916/113033}.

\bibitem{on_the_link_between_binomial_theorem_and_discrete_convolution}
{Kolosov, Petro}.
\newblock {On the link between binomial theorem and discrete convolution}.
\newblock {\em {arXiv preprint arXiv:1603.02468}}, 2016.
\newblock \url{https://arxiv.org/abs/1603.02468}.

\bibitem{unusual_identity_for_odd_powers}
{Kolosov, Petro}.
\newblock {106.37 An unusual identity for odd-powers}.
\newblock {\em {The Mathematical Gazette}}, 106(567):509--513, 2022.
\newblock \url{https://doi.org/10.1017/mag.2022.129}.

\bibitem{history_and_overview_of_polynomial_p}
Petro Kolosov.
\newblock {History and overview of the polynomial P(m,b,x)}, 2024.
\newblock
  \url{https://kolosovpetro.github.io/pdf/HistoryAndOverviewOfPolynomialP.pdf}.

\bibitem{oeis_numerators_of_the_coefficient_a_m_r}
Petro Kolosov.
\newblock {Entry A302971 in The On-Line Encyclopedia of Integer Sequences},
  2018.
\newblock \url{https://oeis.org/A302971}.

\bibitem{oeis_denominators_of_the_coefficient_a_m_r}
Petro Kolosov.
\newblock {Entry A304042 in The On-Line Encyclopedia of Integer Sequences},
  2018.
\newblock \url{https://oeis.org/A304042}.

\bibitem{knuth1993johann}
{Knuth, Donald E.}
\newblock {Johann Faulhaber and sums of powers}.
\newblock {\em {Mathematics of Computation}}, 61(203):277--294, 1993.
\newblock \url{https://arxiv.org/abs/math/9207222}.

\bibitem{kolosov2025unexpected}
Petro Kolosov.
\newblock Unexpected polynomial identity, 2025.
\newblock
  \url{https://kolosovpetro.github.io/pdf/UnexpectedPolynomialIdentity.pdf}.

\bibitem{oeis_coefficients_u_m_l_k_defined_by_polynomial_identity_1}
Petro Kolosov.
\newblock {The coefficients U(m, l, k), m = 1 defined by the polynomial
  identity}, 2018.
\newblock \url{https://oeis.org/A320047}.

\bibitem{oeis_coefficients_u_m_l_k_defined_by_polynomial_identity_2}
Petro Kolosov.
\newblock {The coefficients U(m, l, k), m = 2 defined by the polynomial
  identity}, 2018.
\newblock \url{https://oeis.org/A316349}.

\bibitem{oeis_coefficients_u_m_l_k_defined_by_polynomial_identity_3}
Petro Kolosov.
\newblock {The coefficients U(m, l, k), m = 3 defined by the polynomial
  identity}, 2018.
\newblock \url{https://oeis.org/A316387}.

\bibitem{kolosovpetro_gist}
Petro Kolosov.
\newblock Mathematica programs to generate plots and tables, 2025.
\newblock
  \url{https://gist.github.com/kolosovpetro/2b5c55094c66b8d6a97b9798be9a8dec}.

\end{thebibliography}
